\documentclass[a4paper, reqno, 11pt]{amsart}

\usepackage[english]{babel}
\usepackage{amsmath}
\usepackage{amssymb}
\usepackage{amsthm}
\usepackage{enumerate}
\usepackage{ifthen}
\usepackage{bbm}
\usepackage{color}
\provideboolean{shownotes} 
\setboolean{shownotes}{true}
\usepackage{hyperref}
\usepackage{graphicx}
\usepackage{pgfplots}
\usepackage{tikz,tikz-3dplot} 
\usepackage{mathtools}
\usepackage{comment}
\usepackage{float}
\usepackage{geometry}
\newcommand{\margnote}[1]{
\ifthenelse{\boolean{shownotes}}%
{\marginpar{\raggedright\tiny\texttt{#1}}}%
{}%
}

\newcommand{\hole}[1]{
\ifthenelse{\boolean{shownotes}}%
{\begin{center} \fbox{ \rule {.25cm}{0cm}
\rule[-.1cm]{0cm}{.4cm} \parbox{.85\textwidth}{\begin{center}
\texttt{#1}\end{center}} \rule {.25cm}{0cm}}\end{center}}
{}
}
\usepackage{mathrsfs}
\newtheorem{thm}{Theorem}[section]

\newtheorem{lem}[thm]{Lemma}

\newtheorem{rem}[thm]{Remark}

\newtheorem{defn}[thm]{Definition}

\newcommand{\R}{\mathbb{R}}
\newcommand{\T}{\mathbb{T}}

\newcommand{\dive}{\mathop{\mathrm {div}}}
\newcommand{\curl}{\mathop{\mathrm {curl}}}

\newcommand{\de}{\mathrm{d}}

\numberwithin{equation}{section}

\subjclass[]{}


\begin{document}

\title[On the topology of the magnetic lines of solutions of MHD]{On the topology of the magnetic lines of solutions of the MHD equations}

\author[G. Ciampa]{Gennaro Ciampa}
\address[G.\ Ciampa]{Dipartimento di Matematica ``Federigo Enriques", Universit\`a degli Studi di Milano, Via Cesare Saldini 50, 20133 Milano, Italy.}
\email[]{\href{gciampa@}{gennaro.ciampa@unimi.it}}

\begin{abstract}
We construct examples of smooth periodic solutions to the Magnetohydrodynamic equations in dimension 2 with positive resistivity for which the topology of the magnetic lines changes under the flow. By Alfvén's theorem this is known to be impossible in the ideal case (resistivity = 0). In the resistive case the reconnection of the magnetic lines is known to occur and has deep physical implications, being responsible for many dynamic phenomena in astrophysics. The construction is a simplified proof of \cite{CCL} and in addition we consider the case of the forced system.
\end{abstract}

\maketitle

\section{Introduction}
In these notes we investigate the {\em magnetic reconnection} phenomenon for smooth solutions of Magnetohydrodynamic equations (MHD). These equations describe the behaviour of an electrically conducting incompressible fluid such as plasmas, salt water, liquid metals, etc.. They are obtained by a combination of the incompressible Navier-Stokes equations and the Faraday-Maxwell system via Ohm’s law. The Cauchy problem for the \eqref{eq:mhd} equations is 
\begin{equation}\label{eq:mhd}\tag{MHD}
\begin{cases}
\partial_t u+(u\cdot \nabla)u+\nabla P=\nu\Delta u+(b\cdot \nabla)b,\\
\partial_t b+(u\cdot \nabla)b=(b\cdot \nabla)u+\eta\Delta b,\\
\dive u=\dive b=0,\\
u(0,\cdot)=u_0,\hspace{0.3cm} b(0,\cdot)=b_0,
\end{cases}
\end{equation}
where the unknowns are the magnetic field $b$, the fluid velocity $u$, and the total pressure acting on the fluid $P$. The data of the problem are two divergence-free vector fields $u_0$ and $b_0$, and the parameters $\nu\geq 0$ (the viscosity) and $\eta\geq 0$ (the resistivity). 
For a general review we refer to \cite{Dav, DL72, Priest, ST}. The \eqref{eq:mhd} equations have received considerable attention from mathematicians and the first results on the well-posedness theory date back to \cite{DL72, ST}. The existence of global weak solutions with finite energy and local strong solutions to \eqref{eq:mhd} in two and three dimensions have been proved in~\cite{DL72}. Moreover, for smooth initial data they proved the smoothness and uniqueness of their global weak solutions in the two-dimensional case. On the other hand, in \cite{ST} the authors proved the uniqueness of the local strong solutions in 3D.
In recent years there were different results even for the non-resistive case ($\eta =0$), and local well-posedness results, at an (essentially) sharp level of Sobolev regularity, are now available \cite{FMRR, Feff}. 
Our goal, however, is not to provide a complete list of results for which we refer the reader to the references cited above and to the references contained therein.\\

Our aim is to provide analytical examples of \emph{magnetic reconnection}: it refers to a change in the topology of the magnetic lines, i.e. the integral lines of the magnetic field $b$. We will focus on the two-dimensional periodic case, i.e. solutions defined on the 2D torus $\T^2$. We will say that a solution $(u,b)$ of \eqref{eq:mhd} shows magnetic reconnection if there exist $t_1, t_2$ such that there is no homeomorphism of $\T^2$ mapping the set of the integral lines of $b(t_1,\cdot)$ into integral lines of $b(t_2,\cdot)$. A break in the topology of coherent magnetic structures can release a large amount of energy (that was stored as mechanical or magnetic energy until the magnetic connections evolved coherently), converting it into other forms of energy such as kinetic or thermal energy. This phenomenon is of fundamental importance to astrophysicists being responsible for many dynamic phenomena such as flares, coronal mass ejections, and the solar wind.\\

In the non-resistive case it is known that the integral lines of a sufficiently smooth magnetic field are transported by the fluid ({\em Alfven's theorem}). In particular, the topology of the integral lines of the magnetic field is frozen under the evolution. We briefly explain why: denote by $\Phi_t$ the fluid flow, i.e.
$$
\begin{cases}
\frac{\de}{\de t} \Phi_t (x) = u(t,\Phi_t (x)), \\
\Phi_0 (x) = x,
\end{cases}
$$
and let $b_0$ be a (smooth) magnetic field with $\gamma_0$ being an integral line of $b_0$, i.e.
$$
\frac{\de}{\de s}\gamma_0(s) = b_0(\gamma_0(s)).
$$
If $b(t,x)$ is the unique solution of the PDE
\begin{equation}\label{eq:alfven}
\partial_t b+(u\cdot \nabla) b=  (b\cdot\nabla)u,
\end{equation}
with initial datum $b_0$, it is known that it satisfies the formula
\begin{equation}\label{eq:alfven2}
b(t,\Phi_t(x))=\nabla \Phi_t(x)b_0(x),
\end{equation}
where $\Phi_t$ is the flow of $u$. The formula \eqref{eq:alfven2} defines the so-called {\em pull-back} of $b_0$ by $\Phi_t$. Then, if $\gamma_0$ is an integral curve of $b_0$, we have that
\begin{align*}
\frac{\de}{\de s}\Phi_t(\gamma_0(s))&=\nabla\Phi_t(\gamma_0(s))\frac{\de}{\de s}\gamma_0(s)=\nabla\Phi_t(\gamma_0(s))b_0(\gamma_0(s))=b(t,\Phi_t(\gamma_0(s))),
\end{align*}
meaning that $\gamma_t = \Phi_t \circ \gamma_0$ is an integral line of $b(t,\cdot)$: the integral lines of $b_0$ are transported by the fluid flow. Since $u$ is assumed to be smooth, $\Phi_t:\T^2\to\T^2$ is a diffeomorphism and this implies that, at every time $t>0$, the integral lines of $b_0$ and $b(t,\cdot)$ are diffeomorphic. Thus, there is no reconnection.

Lastly, let us mention that in the non resistive case the {\em topological stability} of the magnetic structures is also related to the conservation of the magnetic helicity, which is a measure of the linkage and twist of the magnetic field lines. 
This feature becomes a very subtle matter at low regularities, intimately related to anomalous dissipation phenomena. We refer to \cite{BBV20,FL19,FLS20}  for some positive and negative results in this direction.\\

In the resistive case ($\eta > 0$) the topology of the magnetic lines it is expected to change under the fluid evolution, even for regular solutions. The heuristic underlying the phenomenon is that magnetic diffusion allows breaking the topological rigidity. It is interesting that we can prove magnetic reconnection at arbitrarily small resistivity, namely for all $\eta >0$, thus even in a very turbulent regime. Although numerical and experimental evidences exist (see \cite{MagnPhys, Priest} and the references therein), no analytical examples of magnetic reconnection are known. 
The main result is the following.
\begin{thm}\label{thm:main}
Given any viscosity and resistivity $\nu, \eta >0$ and any constant $T>0$, there exist a zero-average divergence-free smooth vector field $b_0$ and a unique global smooth solution $(u,b)$ of \eqref{eq:mhd} on $\T^2$ with initial datum $(0, b_0)$, such that the magnetic lines at time $t=0$ and $t=T$ are not topologically equivalent.
\end{thm}

It is known that in two dimensions the reconnection occur only at critical points, see \cite{Priest}. Our proof exploits this feature considering as topological constraint the number of critical points of the magnetic field. Specifically, we will implement a perturbative analysis of some particular solutions of the linearized equation for which one can infer that the reconnection occurs.
In particular, we define the initial magnetic field as
$$
T_{nm} + \delta \tilde T_1,  \quad 0 <  \delta \ll 1,
$$ 
where $T_{nm}$ and $\tilde T_1$ are two Taylor fields (see Section 2 below) with the following properties:
\begin{enumerate}
\item the field $T_{nm}$ has several stagnation points (namely $\sim n^2+m^2\gg 1$, half of them are hyperbolic and half of them elliptic);
\item the field $\tilde T_1$ has exactly four stagnation points, and the topology of the integral lines is completely prescribed (see figure \ref{fig:v1}) and robust ($\tilde T_1$ is structurally stable by Theorem \ref{thm:MW}).
\end{enumerate}

The idea is then to choose the relevant parameters of the construction in such a way that the solution will have at least the same number of critical points of $T_{nm}$ at time $t=0$, while it will be topologically equivalent to $\tilde T_1$ at time $T$. 
This proves that the magnetic reconnection occurred for intermediate times.\\

It is worth noting that we can trivially build examples in the 3D case: just consider a 2D solution that reconnects and define $0$ as the third component. However, it is known that in the three dimensional case the reconnection is not constrained to occur at critical points and several other magnetic field structures may be sites of reconnection, see \cite{Priest}. In this regards, a genuinely 3D argument with the possibility to prescribe rich topological structures is proved in \cite{CCL}, relying upon some deep results about topological richness of Beltrami fields \cite{ELP, Annals, Acta, EPT}. 

Finally, we recall that another very important model for analyzing reconnection phenomena is the Hall-MHD system. The latter has an additional term on the left hand side of the magnetic equation in \eqref{eq:mhd}, namely $\curl(\curl b\wedge b)$, which is called the {\em Hall term}. We refer to \cite{ADFL, CDL, DS, HG} and reference therein for an introduction to these equations. 
The Hall term alone cannot change the topology of the magnetic field lines, however it interacts with the magnetic viscosity accelerating the reconnection process, see \cite{HG}. It is interesting to note that the Hall term is identically zero if $b$ is a Beltrami field. Thus, since a small data theory is available for global smooth solutions of the Hall-MHD \cite{CDL}, we expect that the perturbative argument of the 3D proof in \cite{CCL} can be adapted to the case of the Hall-MHD system as well.
\\

Finally, we will show how adding a force to the system \eqref{eq:mhd} can help construct an (explicit) solution that exhibits magnetic reconnection. The theorem is the following.

\begin{thm}\label{thm:main2}
Given any viscosity and resistivity $\nu, \eta >0$ and any constant $T>0$, there exists a smooth zero-average and divergence-free vector field $b_0$ and a couple $(f_1,f_2)$ of smooth vector fields with zero average, such that:
\begin{itemize}
\item there exists a zero-average unique global smooth solution $(u,b)$ of the \eqref{eq:mhd} system with forces $(f_1,f_2)$ and initial datum $(0, b_0)$;
\item the magnetic lines at time $t=0$ and $t= T$ are not topologically equivalent.
\end{itemize}
\end{thm}
Differently from the proof of Theorem \ref{thm:main}, the proof of Theorem \ref{thm:main2} is not based on a perturbative argument and we will explicitly construct the forces and the solution. 

\section{Taylor fields}\label{Sec:Taylor}
In this section, after recalling some basic definitions, we introduce the building blocks of our reconnection result. We say that $v:\T^2 \to \R^2$ is a {\em Hamiltonian} vector field if it can be expressed as the orthogonal gradient of a scalar function $\psi$, i.e.
$
v=\nabla^\perp \psi:=(\partial_{x_2}\psi,-\partial_{x_1}\psi).
$
Hamiltonian vector fields are by definition divergence-free and 
by the Helmoltz decomposition an incompressible vector field on $\T^2$ is either Hamiltonian or constant. Since we are interested in zero-average vector fields for us there will be no difference between being incompressible and Hamiltonian.
We recall that a singular point $x_0$ of a vector field $v\in C^1(\T^2)$ is said to be {\em non-degenerate} if $\nabla v(x_0)$ is an invertible matrix. A non-degenerate singular point of a divergence-free vector field must be either a saddle or a center. We now give the definition of structural stability.
\begin{defn}\label{def:struc-stab}
A vector field $v$ on $\T^d$ is structurally stable if there is a neighborhood $\mathcal{U}$ of $v$ in $C^1(\T^d)$ such that whenever $v'\in \mathcal{U}$ there is a homeomorphism of $\T^d$ onto itself transforming trajectories of $v$ onto trajectories of $v'$.
\end{defn}

As a consequence of the classical Peixoto's Theorem, divergence-free vector fields are not topologically stable under arbitrary perturbations. However, the following characterization of structurally stability holds for Hamiltonian vector fields on the two-dimensional torus.

\begin{thm}[Ma, Wang \cite{MW}]\label{thm:MW}
A divergence-free Hamiltonian vector field $v\in C^r(\T^2)$ with $r\geq 1$ is structurally stable under Hamiltonian vector field perturbations if and only if
\begin{itemize}
\item all singular points of $v$ are not degenerate;
\item all saddle connections of $v$ are self saddle connections.
\end{itemize}
\end{thm}

\begin{figure}[H]
\centering
\includegraphics[width=0.5\textwidth]{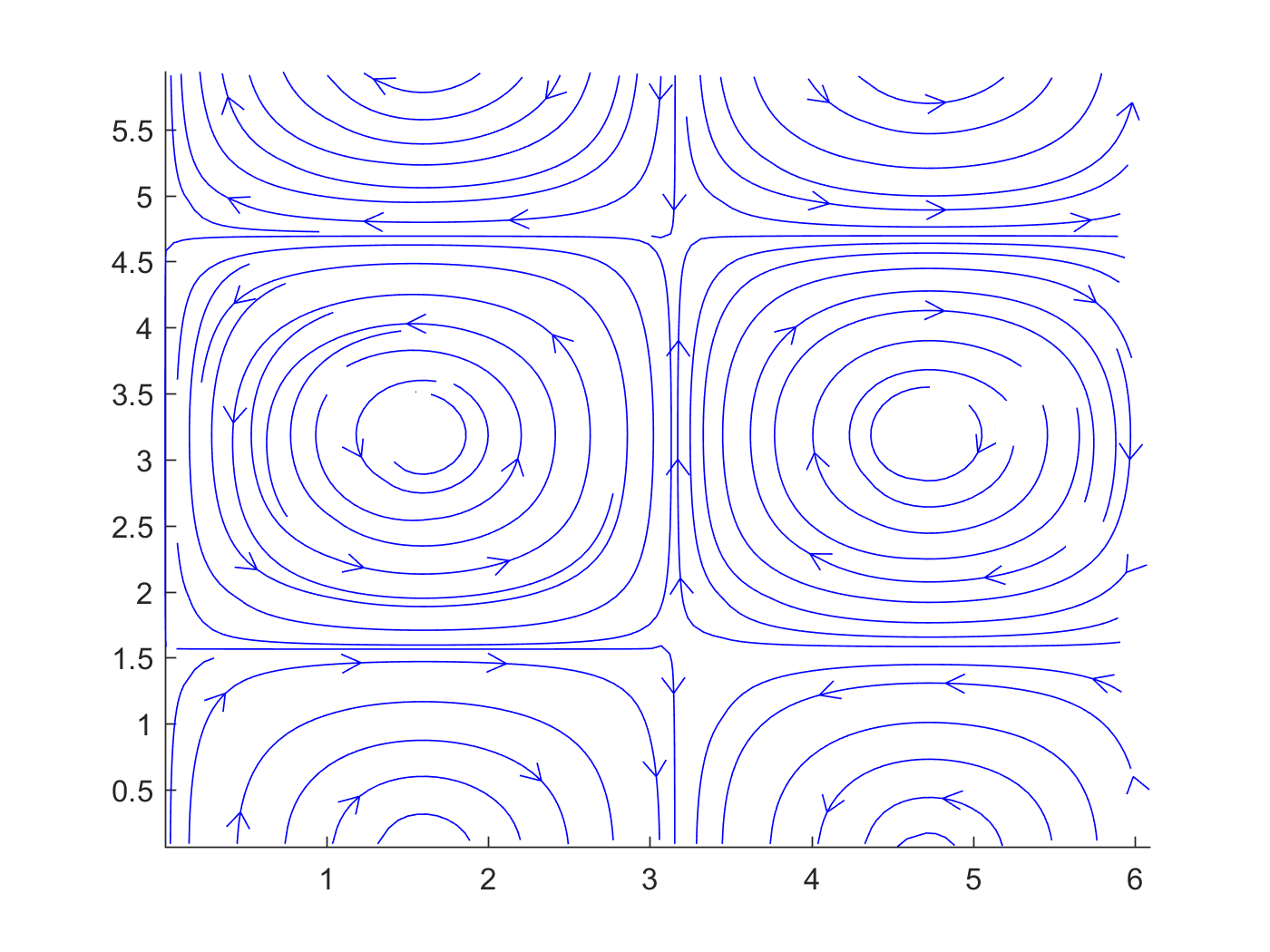}
\caption{Integral lines of $T_{11}$.}
\label{fig:v22}
\end{figure}

Taylor fields on $\T^2$ are divergence-free eigenfunctions of the Laplacian, i.e. solutions of
$$
\begin{cases}\label{def:eigenvalue-taylor}
-\Delta T_M=M\, T_M,\\
T_M=\nabla^{\perp} h_M,
\end{cases}
$$
with $M \neq  0$. They are Hamiltonian with stream function $h_M$, where $\{ h_k \}_k$ form a orthonormal basis of $L^2(\T^2)$. We will focus on some particular choices:
\begin{equation}\label{def:Tnm}
T_{nm} = (m\sin nx \sin my, n\cos nx \cos my),
\end{equation}
with eigenvalue $n^2+m^2$, and
\begin{equation}\label{def:T1}
\tilde T_1=\left(\sin y,\frac{1}{2}\sin x\right),
\end{equation}
with eigenvalue $1$. It is important to note that the vector fields $T_{nm},\tilde T_1$ are also stationary solutions of the 2D Euler equations for a suitable choice of the pressure. Moreover, the field $T_{nm}$ is not structurally stable in the sense of Definition \ref{def:struc-stab}. The instability follows from the presence of saddle connections in the phase diagram of $T_{nm}$, see Figure \ref{fig:v22}.

\begin{figure}[H]
\centering
\includegraphics[width=0.5\textwidth]{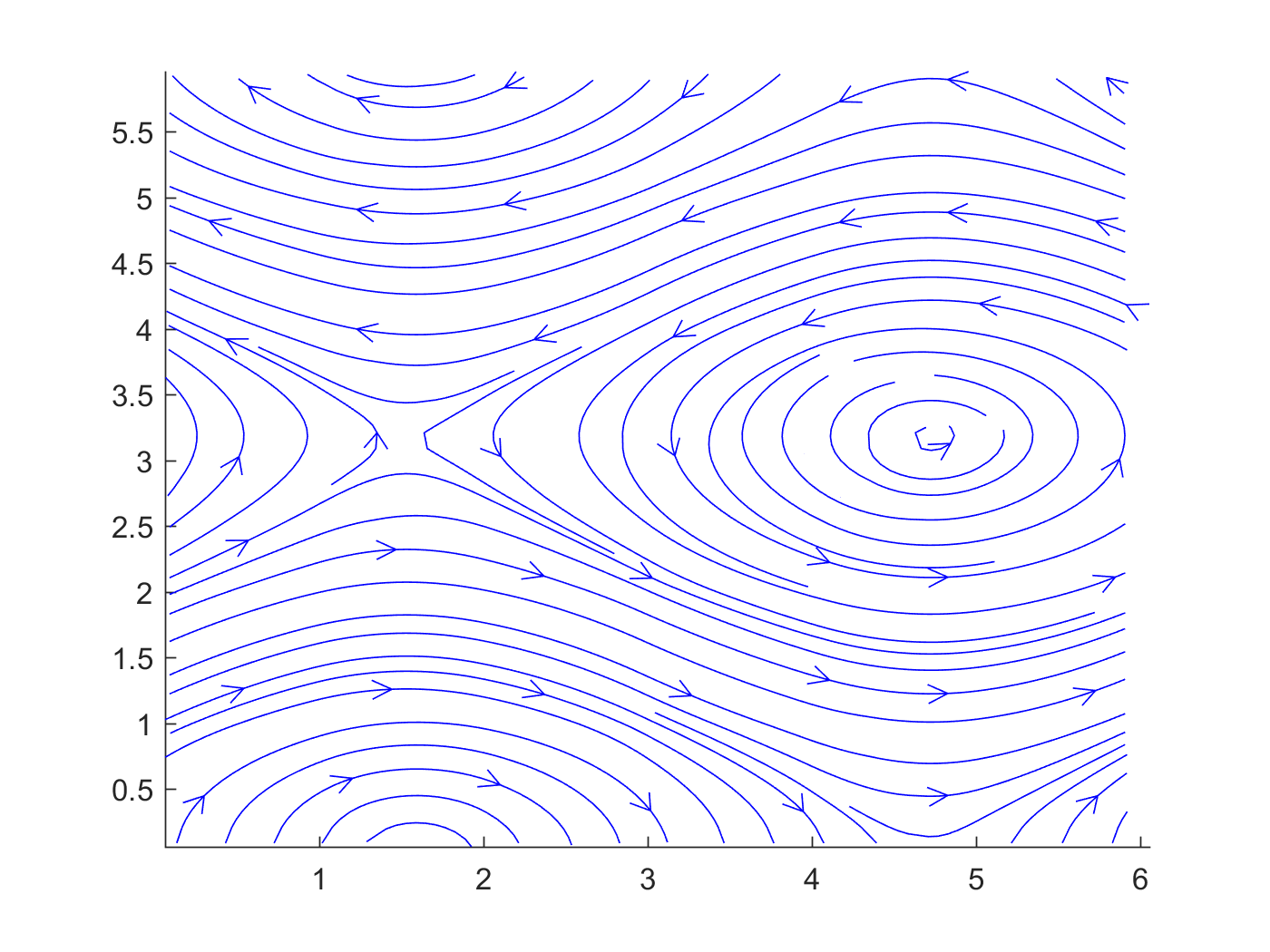}
\caption{Integral lines of $\tilde T_{1}$.}
\label{fig:v1}
\end{figure}

Instead, the vector field $\tilde T_1$ is structurally stable as there are no saddle connections in its phase diagram except self-connections. The field $\tilde T_1$ has indeed two saddle points which are not connected by any integral line, see Figure \ref{fig:v1}.

We conclude this subsection with the lemma below. It investigates the stability of the critical points of the Taylor fields and in particular it shows that the number of critical points does not decrease if the perturbation is small enough (with a quantitative bound), see \cite{CCL}. 

\begin{lem}\label{Lemma:CriticalPoints}
Let $T_{nm}$ be the Taylor field with eigenvalue $N^2=n^2+m^2$ defined above and let $ W \in C^1$. For all $N$ sufficiently large (depending only on $\| W \|_{C^1}$), there exists $\delta_0=\delta_0(N)$ such that the vector field
$
\tilde{V}(x):=T_{nm}(x)+\delta  W(x),
$
has at least $8nm$ regular critical points for every $|\delta|<\delta_0(N)$. We may choose $\delta_0(N) =  N^{-L}$, where~$L$ is a fixed large integer. 
\end{lem}

\subsection{Stability estimates for the MHD equations}
Here we state a stability result for the \eqref{eq:mhd}. We consider two solutions $(u,b)$ and $(w,m)$ of \eqref{eq:mhd} with initial datum, respectively, $(u_0,b_0)$ and $(w_0,m_0)$. By defining the differences $v=u-w$, $h=b-m$, we obtain that the couple $(v,h)$ solves the system 
\begin{equation}\label{eq:differenceTris}
\begin{cases}
\partial_t v+\dive\left(v\otimes v+2v\otimes w\right)+\nabla P_{v,h}=\nu\Delta v+\dive\left(h\otimes h+2h\otimes m\right),\\
\partial_t h+(v\cdot\nabla )h+(w\cdot \nabla)h+(v\cdot\nabla)m=\eta\Delta h+(h\cdot\nabla)v+(h\cdot\nabla)w+(m\cdot\nabla)v, \\
\dive v=\dive h=0,\\
v_0=u_0-w_0,\hspace{0.3cm}h_0=b_0-m_0.
\end{cases}
\end{equation}
Define 
$
\tilde \Gamma :=  1 +
\|w_0\|_{L^2}^2+\|m_0\|_{L^2}^2 
+ \|v_0\|_{L^2}^2+\|h_0\|_{L^2}^2
$
and take $\sigma < \min ( \nu,\eta )$. The following theorem holds, see \cite{CCL}.

\begin{thm}\label{thm:boundhrStab}
Let $v_0,h_0\in H^r(\T^2)$ be two divergence-free vector fields with zero mean. Let $(v,h,P_{v,h})$ be the unique solution of \eqref{eq:differenceTris} with initial datum $(v_0,h_0)$, where $(w,m)$ is a solution of \eqref{eq:mhd}. Assume that 
$$
\| w_0 \|_{H^m} + \| m_0 \|_{H^m}  \leq C N^m, \qquad \| v_0 \|_{H^m} + \| h_0 \|_{H^m} \leq C \delta, \qquad m=0, \ldots, r,
$$  
for some $N > 1$. Then
\begin{equation}\label{est:HrStab}
\|h(t,\cdot)\|_{H^r}^2 + \|v(t,\cdot)\|_{H^r}^2
\leq  
\delta^2 N^{2r} e^{-2\sigma t} e^{\frac{C \tilde\Gamma}{\sigma^2}},
\end{equation}
where the implicit constants $C$ depend on $r, \sigma$.
\end{thm}

\section{Proof of the reconnection}
In this section we prove our first main theorem.
\begin{proof}[Theorem \ref{thm:main}]
We consider as initial datum 
$$
(u_0,b_0):=\left(0, N^{-1}T_{nm}+  \delta \tilde T_1\right)
$$
with $n^2 + m^2 \gg 1$ and $\delta \ll 1$ to be chosen later. Consider the rescaled datum
$
N b_0=T_{nm}+ N \delta \tilde T_1,
$
and take
$
 \delta < \frac{c}{N^{L+1}},
$
where $c$ is a suitable small constant, we have that $\delta \ll \delta_{0}(N)$ from Lemma \ref{Lemma:CriticalPoints}, thus the vector field $Nb_0$, and so $b_0$, has at least $8nm$ regular critical points. With this choice,  at time $t=0$ the solution starts from a configuration which is not topological equivalent to $\tilde T_1$.\\

Now we will show that at time $t=T$ (a rescaled version of) the solution $b(T,\cdot)$ and the vector field $\tilde T_1$ are close in the $C^1$-norm: the structural stability of the latter will imply that $b(T,\cdot)$ is topologically equivalent to $\tilde T_1$ and the magnetic reconnection happened between $t = 0$ and $t = T > 0$. To this end, we will make use of a perturbative argument with respect to a reference (given) solution that we now construct.\\

Let us consider the initial datum $(w_0,m_0)=\left(0,N^{-1}T_{nm}\right)$ for \eqref{eq:mhd}. We recall that the Taylor field $T_{nm}$ is a solution of stationary Euler 
with a suitable choice of pressure that we denote by $P_{T_{nm}}$.
Moreover, by definition $N^2 = n^2 + m^2$ and
$
\Delta T_{nm} = -N^2 T_{nm}.
$
Then, it is easy to check that the unique solution of \eqref{eq:mhd} with initial datum $(w_0,m_0)$ is given by
\begin{equation}\label{def:wm}
(w,m,p_w) = \left(0, N^{-1}e^{-\eta N^2 t} T_{nm},N^{-2}e^{-2\eta N^2 t} P_{T_{nm}} \right).
\end{equation}

We consider the behavior of the fluid at time $t=T$. In order to deal with the non-linear contribution of the equations, we define $u=v+w,b=h+m$ where $(v,h)$ satisfies the difference equation \eqref{eq:differenceTris}. By Duhamel's formula, we know that the solution $b$ satisfies the identity
\begin{equation}
b(t,\cdot)=e^{t\Delta}b_0+D(t,\cdot),
\end{equation}
where 
$
D(t,\cdot) := L_h(t,\cdot) + L_b(t,\cdot),
$
and
\begin{equation}\label{Def:Lh}
L_h(t,\cdot):=\int_0^t e^{\eta(t-s)\Delta}\dive \big(h(s)\otimes v(s)-v(s)\otimes h(s)\big)\,\de s,
\end{equation}
\begin{equation}\label{Def:Lb}
L_m(t,\cdot):=\int_0^t e^{\eta(t-s)\Delta}\dive \big(m(s)\otimes v(s)-v(s)\otimes m(s)\big)\,\de s.
\end{equation}
We have the following Lemma, see \cite{CCL}.
\begin{lem}\label{lem:duhamel}
Let $(v,h)$ be the unique solution of \eqref{eq:differenceTris} with initial datum $(0, \delta \tilde T_1)$, whith $(w,m)$ defined in \eqref{def:wm}. Then, the functions $L_h$, $L_m$ defined in \eqref{Def:Lh} and \eqref{Def:Lb} satisfy the following estimates
$$
\|L_h(t,\cdot)\|_{H^r}\leq C \delta^2 N^{r+3} e^{-\sigma s},\,\,
\|L_m(t,\cdot)\|_{H^r}\leq C\delta N^{-2}+ C\delta  N^{r+1} e^{- \eta N^2 \frac{t}{2}}.
$$
\end{lem}
The proof of the previous Lemma follows from an application of Theorem \ref{thm:boundhrStab}: this is why we consider the factor $N^{-1}$ in the definition of $b_0$, but clearly the proof can be adapted by slightly changing the exponents of the parameter $N$.
Rescale the magnetic field $b(T,\cdot)$ by the factor $\delta^{-1} e^{\eta T}$ and denote it by $\tilde b(T,\cdot)$. More explicitely
$$
\tilde b(T,\cdot) =  \tilde T_1 +(\delta N)^{-1}e^{-\eta (N^2 - 1) T} T_{nm}  + \delta^{-1} e^{\eta T} D(T,\cdot).
$$
Our goal is to choose $N$ so large such that  
$
\left\| \tilde b(T,\cdot) -  \tilde T_1  \right\|_{H^r} \ll 1,
$
and we will use the structural stability of $\tilde T_1$ under Hamiltonian perturbations and Sobolev embeddings ($r=3$) to infer that the set of the integral lines of $\tilde b(T,\cdot)$, and thus of $b(T,\cdot)$, is diffeomorphic to that of $\tilde T_1$. For some sufficiently large $L$, we choose  
$$
\delta =  e^{- \eta T} N^{-(L+1)}, \quad L \geq r+3.
$$ 
This is compatible with Lemma \ref{Lemma:CriticalPoints} and then combining it with Lemma \ref{lem:duhamel}, we obtain
$$
\| \delta^{-1} e^{\eta T} D(t,\cdot) \|_{H^r} \leq CN^{-1} + CN^{-2} e^{\eta T} + C N^{r+1} e^{- \eta N^2 \frac{T}{2}},
$$
which is small if $N$ is sufficiently large, depending on $T$ and $\eta$. With this choice of $\delta$, up to taking $N$ even larger, we can also have that 
$$
\| \delta^{-1} e^{-\eta (N^2-1) T} T_{nm} \|_{H^r} \leq 
C N^{L+1+r}e^{-\eta N^2 T} \ll 1.
$$
This shows that at time $T$ the solution is topological equivalent to $\tilde T_1$ and then the reconnection of the magnetic lines took place.

\end{proof}

\section{The forced system}
In this section we prove our second main result, namely Theorem \ref{thm:main2}.

\begin{proof}[Theorem \ref{thm:main2}]
As already mentioned in the introduction, we are going to explicitly construct $f_1$, $f_2$, $u$, $b$. We define $b_0= T_{nm}$ and $f_2=T_{N_2}$ where $T_{N_2}$ is a Taylor field of the same form of \eqref{def:Tnm} but with eigenvalue such that $N_2<<N$. In particular, the number of critical points of $T_{nm}$ is much greater than those of $T_{N_2}$. Moreover, notice that $f_2$ has zero mean. The idea now is to choose $f_1$ such that the solution is of the form $(0,b(t,x))$. By Duhamel's formula, the solution of the forced heat equation
$$
\begin{cases}
\partial_t b=\eta\Delta b+ f_2\\
b(0,\cdot)=T_{nm},
\end{cases}
$$
is given by 
$$
b(t,x)= e^{-\eta N^2 t} T_{nm}(x)+\left(\frac{1-e^{-N_2^2\eta t}}{\eta N_2^2}\right) T_{N_2}(x),
$$
where we used that Taylor fields are eigenfunctions of the Laplacian.
We compute $(b\cdot\nabla)b$: since $T_{nm}$ and $T_{N_2}$ are stationary solutions of 2D Euler, we obtain that
\begin{align*}
(b\cdot\nabla)b
&=\underbrace{ e^{-2\eta N^2 t} \nabla P_{nm}+\left(\frac{1-e^{-\eta N_2^2 t}}{\eta N_2^2}\right)^2\nabla P_2}_{(*)}\\
&+\underbrace{ e^{-\eta N^2 t}\left(\frac{1-e^{-\eta t}}{\eta}\right)\left((T_{nm}\cdot\nabla) T_{N_2}+( T_{N_2}\cdot\nabla)T_{nm}\right).}_{(**)}
\end{align*}
Then, if we define $\nabla P=(*)$ and $f_1=-(**)$, we have that the triple $(0,b,P)$ solves the forced \eqref{eq:mhd} with $f_1$ having zero mean. Now we are in position to prove the reconnection: by defining the rescaled magnetic field $\tilde b=\eta N^2_2 b$, we have that
\begin{equation}
\|\tilde b(T,\cdot)- T_{N_2}\|_{H^r}\leq C\eta N^r N_2^2 e^{-\eta N^2 T}+CN_2^r e^{-\eta N_2^2T},
\end{equation}
and then the conclusion follows from Lemma \ref{Lemma:CriticalPoints} by choosing  $N_2$ first and then $N>>N_2$.
\end{proof}

\begin{rem}
It is clear from the previous proof that the choice of $T_{N_2}$ is not limited to Taylor fields of the form \eqref{def:Tnm}. In fact, other Taylor fields can be considered as long as the condition on the eigenvalues $N_2<<N$ is satisfied.
\end{rem}

\begin{rem}
Note that we can also give a different proof provided that we are free to choose an arbitrarily large time $T$. Define $f_1=\tilde T_1$ where $\tilde T_1$ is defined in \eqref{def:T1}. With the same computations as above we obtain that, if $\tilde b=\eta b$,
$$
\|\tilde b(T,\cdot)-\tilde T_1\|_{H^r}\leq \eta ( N^r e^{-\eta N^2 T}+e^{-\eta T}),
$$
which can be made as small as we want for $T$ big enough. Then, the proof follows from the structurally stability of $\tilde T_1$.
\end{rem}

\section*{Acknowledgements}
The author is supported by the ERC STARTING GRANT 2021 ``Hamiltonian Dynamics, Normal Forms and Water Waves" (HamDyWWa), Project Number: 101039762. Views and opinions expressed are however those of the authors only and do not necessarily reflect those of the European Union or the European Research Council. Neither the European Union nor the granting authority can be held responsible for them. This research is partially funded by AEI through the project PID2021-123034NB-I00.

\end{document}